\documentclass[10pt]{amsart}
\usepackage{amssymb}
\usepackage{epsfig}
\usepackage{amsfonts}
\usepackage{amscd}

\newcommand{\z}{\mathbb{Z}}

\newcommand{\cx}{\mathbb{C}}

\newtheorem{thm}{\textbf{Theorem }}[section]
\newtheorem{prop}[thm]{\textbf{Proposition }}%[section]
\newtheorem{cor}[thm]{\textbf{Corollary}}%[section]
\newtheorem{lemma}[thm]{\textbf{Lemma}}%[section]

\newcommand{\eqqh}{{QH}^{*}_{T}(Gr(p,m))}
\newcommand{\eqqhx}{{QH}^{*}_{T}(X)}

\numberwithin{equation}{section}

\begin{document}
\title[Giambelli formulae for equivariant quantum cohomology]{Giambelli formulae for the equivariant quantum cohomology of the Grassmannian}
\author{Leonardo Constantin Mihalcea} \thanks{MSC: 14N35 (primary); 05E05;14F43 (secondary).}\address{525 E. University,
Dept. of Mathematics, University of Michigan, Ann Arbor, MI
48109}\email{lmihalce@umich.edu}
\date{\today} \maketitle
\begin{abstract} We find presentations by generators and relations for the equivariant quantum cohomology of the Grassmannian. For these presentations, we also find determinantal formulae for the equivariant quantum Schubert classes. To prove this, we use the theory of factorial Schur functions and a characterization of the equivariant quantum cohomology ring. \end{abstract}
\section{Introduction}

Let $X$ denote the Grassmannian $Gr(p,m)$ of subspaces of
dimension $p$ in $\cx^m$. One of the fundamental problems in the study of the equivariant quantum cohomology algebra of $X$ is to compute its structure constants, which are the 3-point, genus 0, equivariant Gromov-Witten invariants. The goal of this paper is to give a method for such a computation. Concretely, we realize the equivariant quantum cohomology as a ring given by generators and relations and we find polynomial representatives (i.e. Giambelli formulae) for the equivariant quantum Schubert classes (which form a module-basis for this ring).\begin{footnote}{It is a standard fact that given a ring $R/I$, where $R$ is a polynomial ring and $I$ an ideal, together with some elements in $R$ which determine a module-basis for $R/I$, the structure constants for this basis can be computed using e.g. Groebner basis methods.}\end{footnote} These polynomials will be given by certain determinants which appear in
the Jacobi-Trudi formulae for the factorial Schur functions (see \S 2 below for details). 

Since the equivariant quantum cohomology ring specializes to both
quantum and equivariant cohomology rings, we also obtain, as
corollaries, determinantal formulae for the Schubert classes in
the quantum and equivariant cohomology. In fact, in the quantum
case, we recover Bertram's quantum Giambelli formula \cite{Be}. In
the case of equivariant cohomology, we show that the factorial
Schur functions represent the equivariant Schubert classes. The
latter result, although not explicitly stated in the literature,
seems to have been known before (cf. Remark 2 in \S 5).

We recall next some of the basic facts about the equivariant quantum cohomology and fix the notation. The torus $T \simeq (\cx^*)^m$ acts on the
Grassmannian $X$ by the action induced from the $GL(m)-$action.
The $T-$equivariant cohomology of a point, denoted $\Lambda$, is
the polynomial ring $\z[T_1,...,T_m]$ in the equivariant
parameters $T_i$, graded by $\deg T_i = 1$ (see \S 5.1 for a geometric interpretation of $T_i$). Let $q$ be an
indeterminate of degree $m$. The $T-$equivariant quantum
cohomology of $X$, denoted $QH^*_T(X)$, is a graded, commutative,
$\Lambda[q]-$algebra with a $\Lambda[q]-$basis
$\{\sigma_\lambda\}$ indexed by partitions
$\lambda=(\lambda_1,...,\lambda_p)$ included in the $p \times
(m-p)$ rectangle (i.e. $\lambda_1,...,\lambda_p$ are
integers such that $m-p \geqslant \lambda_1 \geqslant ...
\geqslant \lambda_p \geqslant 0$). This basis is determined by the
Schubert varieties in $X$, defined with respect to the standard
flag, and the classes $\sigma_\lambda$ will be called (equivariant quantum) {\it
Schubert classes}. More details are given in \S \ref{sec
eqq} and especially in \cite{Mi1}. The equivariant quantum
multiplication is denoted by $\circ$ and it is determined by the
3-pointed, genus $0$, equivariant Gromov-Witten invariants
$c_{\lambda,\mu}^{\nu,d}$. In this paper we refer to these
coefficients as the equivariant quantum Littlewood-Richardson
coefficients, abbreviated EQLR. They have been introduced by
Givental-Kim in \cite{GK} (see also \cite{Kim2,G,Kim3}), together
with the equivariant quantum cohomology. Then, by definition,
\[ \sigma_\lambda \circ \sigma_\mu= \sum_{d \geqslant 0} \sum_\nu
q^d c_{\lambda,\mu}^{\nu,d} \sigma_\nu. \] The EQLR coefficients
$c_{\lambda,\mu}^{\nu,d}$ are homogeneous polynomials in $\Lambda$
of degree $|\lambda|+|\mu|-|\nu|-md$, where $|\alpha|$ denotes the
sum of all parts of the partition $\alpha$. They are equal to the
structure constants of equivariant cohomology if $d=0$ and to
those of quantum cohomology (i.e. to the ordinary Gromov-Witten invariants) if $|\lambda|+|\mu|-|\nu|-md=0$. Equivalently, the quotient of equivariant quantum cohomology ring
by the ideal generated by the equivariant parameters $T_i$ 
yields
the quantum cohomology ring of $X$ (a $\z[q]-$algebra), while the
quotient by the ideal generated by $q$ yields the $T-$equivariant
cohomology of $X$ (a $\Lambda-$algebra). More about the
EQLR coefficients, including a certain positivity, which
generalizes the positivity enjoyed by the
equivariant coefficients, can be found in \cite{Mi1,Mi2}.

%
%One of the fundamental problems in the study of the equivariant quantum cohomology
%algebra is to find formulae or algorithms to compute its structure constants, the EQLR coefficients $c_{\lambda,\mu}^{\nu,d}$. Such an algorithm can be found in \cite{Mi1}, using a certain recurrence formula. In this paper we follow a more traditional path, in Schubert calculus, by exhibiting a presentation 

%about the  

%More about the
%EQLR coefficients, including a certain positivity, which
%generalizes the positivity enjoyed by the corresponding
%equivariant coefficients, can be found in \cite{Mi1,Mi2}.

\subsection{Statement of the results}\label{subsec statement of res} As we have noted before, the equivariant quantum Giambelli formula which we obtain is a ``factorial"
generalization of Bertram's quantum Giambelli formula \cite{Be}, and, as in that case, it doesn't involve the quantum parameter $q$.
It is closely related to a certain generalization of ordinary Schur
functions, called {\em factorial Schur functions}. These are
polynomials $s_{\lambda}(x;t)$ in two sets of variables:
$x=(x_1,...,x_p)$ and $t=(t_i)_{i \in \z}$. They play a fundamental role in the study of central elements in the universal enveloping algebra of $\frak{gl}(n)$ (\cite{O,OO,Mo}). One of their definitions is via a ``factorial Jacobi-Trudi" determinant, and this
determinant will represent the equivariant quantum Schubert
classes. The basic properties of the factorial Schur functions are
given in \S \ref{fSchurs}.

Denote by $h_i(x;t)$ (respectively by $e_j(x;t)$) the complete
homogeneous (respectively, elementary) factorial Schur functions.  They are equal to
$s_{(i)}(x;t)$ (respectively $s_{(1)^i}(x;t)$), where $(i)$ (resp. $(1)^i$) denotes the partition $(i,0,...,0)$ (resp. $(1,...,1,0,...,0)$, with $i$ $1$'s). Let $t=(t_i)_{i \in \z}$ be the
sequence defined by
$$t_i = \left \{
\begin{array}{ll}
T_{m-i+1} & \textrm{if } 1 \leqslant i \leqslant m \\
0 & \textrm { otherwise}
\end{array} \right.  $$ where $T_i$ is an equivariant parameter. The next definitions are inspired from those used in the theory of factorial Schur functions. For an integer $s$ we define the shifted sequence $\tau^st$ to be
the sequence whose $i-$th term $(\tau^st)_i$ is equal to
$t_{s+i}$. Let $h_1,...,h_{m-p}$ and $e_1,...,e_p$ denote two sets
of indeterminates (these correspond to the complete homogeneous,
respectively, elementary, factorial Schur functions). Often one considers shifts
$s_\lambda(x|\tau^st)$ of $s_\lambda(x|t)$ by shifting the sequence $(t_i)$.
Corresponding to these shifts we define the shifted indeterminates
$\tau^{-s}h_i$, respectively $\tau^s e_j$, where $s$ is a
nonnegative integer, as elements of $\Lambda[h_1,...,h_{m-p}]$,
respectively $\Lambda[e_1,...,e_p]$. 

The definition of $\tau^{-s}h_j$ is given inductively as a function of
$\tau^{-s+1}h_j$ and $\tau^{-s+1}h_{j-1}$, and it is modelled
on an equation which relates $h_j(x|\tau^{-s}t)$ to
$h_j(x|\tau^{-s+1}t)$ (see eq. (\ref{htau}) below, with $a:=\tau^{-s+1}t$).
Concretely, $\tau^0h_j=h_j$,
$\tau^{-1}h_j=h_j+(t_{j-1+p}-t_0)h_{j-1}$ and, in general,
\begin{equation}\label{eq def htau introx}
\tau^{-s}h_j= \tau^{-s+1}h_j + (t_{j+p-s}-t_{1-s})
\tau^{-s+1}h_{j-1}.
\end{equation}
Similarly, the definition of $\tau^s e_i$ is modelled on an
equation which relates $e_{j+1}(x|\tau^{s}t)$ to
$e_{j+1}(x|\tau^{s-1}t)$ (see eq. (\ref{etau}) below, with $a:=\tau^{s-1}t$), and it is given
by
\begin{equation}\label{eq def etau introx}\tau^s e_i = \tau^{s-1}e_i + (t_s -
t_{p-i+s+1}) \tau^{s-1}e_{i-1} \end{equation} with $\tau^0 e_i =
e_i$. By convention, $h_0=e_0 =1$, $h_j=0$ if $j<0$ or $j>m-p$, and
$e_i = 0$ if $i < 0$ or $i>p$.

For $\lambda$ a partition in the $p \times (m-p)$ rectangle define $s_\lambda \in \Lambda[h_1,...,h_{m-p}]$, respectively $\widetilde{s}_\lambda \in \Lambda[e_1,...,e_p]$, analogously to the definition of the factorial Schur function $s_\lambda(x|t)$, via the factorial Jacobi-Trudi determinants (cf.  (\ref{fjt}) below):
\begin{equation}\label{sh introx} s_\lambda=
\det(\tau^{1-j}h_{\lambda_i+j-i})_{1 \leqslant i,j \leqslant p},
\end{equation} 
\begin{equation}\label{se introx} \widetilde{s}_\lambda=\det(\tau^{j-1}e_{\lambda_i'+j-i})_{1
\leqslant i,j \leqslant m-p}. \end{equation} Here $\lambda'=(\lambda_1',...,\lambda_{m-p}')$ denotes the partition
conjugate to $\lambda$, i.e. the partition in the $(m-p) \times p$
rectangle whose $i-$th part is equal to the length of the $i-$th
column of the Young diagram of $\lambda$.

By $H_k$, for $m-p < k \leqslant m$, respectively $E_k$, for $p < k \leqslant m$, we denote the determinants from (\ref{se introx}), respectively (\ref{sh introx}), above, corresponding to partitions $(k)$, respectively $(1)^k$, for the appropriate $k$:
\begin{equation}\label{eq def Hk introx} H_k = \det (\tau^{j-1}
e_{1+j-i})_{1 \leqslant i,j \leqslant k}, \end{equation}

\begin{equation}\label{eq def Ek introx}
E_k = \det (\tau^{1-j}h_{1+j-i})_{1 \leqslant i,j \leqslant k}.
\end{equation} 

With this notation, we present the main result of this paper:
%It is known that $s_{(1)^k}$

\begin{thm}\label{Main Theorem introx}
\noindent (a) There is a canonical isomorphism of
$\Lambda[q]-$algebras
\[\Lambda[q][h_1,...,h_{m-p}]/ \langle E_{p+1},..., E_{m-1}, E_m
+(-1)^{m-p}q \rangle \longrightarrow \eqqhx, \] sending $h_i$ to
$\sigma_{(i)}$. More generally, the image of $s_\lambda$ is the
Schubert class $\sigma_\lambda$.

\noindent (b)(\textrm{Dual version}) There is a canonical
isomorphism of $\Lambda[q]-$algebras
\[\Lambda[q][e_1,...,e_{p}]/ \langle H_{m-p+1},..., H_{m-1}, H_m
+(-1)^{p}q \rangle \longrightarrow \eqqhx, \] sending $e_j$ to
$\sigma_{(1)^j}$ and $\widetilde{s}_\lambda$ to the Schubert class
$\sigma_\lambda$.
\end{thm}

Presentations by generators and relations for the equivariant quantum cohomology of the type A
partial flag manifolds have also been obtained by A. Givental and B. Kim in \cite{GK,Kim1}, 
and, in general, for flag manifolds $X=G/B$ ($G$ a connected, semisimple,
complex Lie group and $B$ a Borel subgroup) by B. Kim in \cite{Kim3}.  In
the case of the Grassmannian, their presentation is given as
\[ \Lambda[a_1,...,a_p, b_1,....,b_{m-p}]/I\] where $I$ is the ideal generated by
$\sum_{i+j=k}a_ib_j=e_k(T_1,...,T_m)$  and $a_pb_{m-p}=e_m(T_1,...,T_m)+(-1)^{p}q$;
the sum is over integers $i,j$
such that $0\leqslant i \leqslant p$ and $0 \leqslant j
\leqslant m-p$; $k$ varies between $1$ and $m-1$ and
$a_0=b_0=1$; $e_k$ denotes the elementary symmetric function.
We also note the results from \cite{Mi1}, where a recursive relation, derived from a multiplication rule with the class $\sigma_{(1)}$ (a Pieri-Chevalley rule), gives another method to compute the EQLR coefficients.

\subsection{Idea of proof} The proof of the theorem uses the theory of factorial Schur functions
and a characterization of the
equivariant quantum cohomology (\cite{Mi1}). We prove first the ``dual version" of the statement. For that we show an equivalent result, where in all the formulae
$\tau^se_j$ is replaced by $e_j(x|\tau^st)$ (here 
$x=(x_1,...,x_p)$ and $t=(t_i)$ is the sequence defined in the
beginning of \S \ref{subsec statement of res}). A key role in this ``translation" is played by a factorial version of the Jacobi-Trudi formula (cf. \cite{Mc1}, pag. 56). Then we prove that the images of the polynomials $s_\lambda(x|t)$, for $\lambda$ included in the $p \times (m-p)$ rectangle, form a $\Lambda[q]-$basis in the claimed presentation. A special
multiplication formula, due to Molev and Sagan (\cite{MS}, see
also \cite{KT1}), computes the product $s_\lambda(x|t)
s_{(1)}(x|t)$ as a sum of $s_\mu(x|t)$, but with $\mu$ having possibly a part larger than $m-p$. Using again the factorial Jacobi-Trudi formula, we prove that, modulo the relations ideal, this multiplication is precisely the equivariant quantum Pieri-Chevalley rule (see \cite{Mi1}). But this rule determines
completely $\eqqhx$, and so the ``dual" statement is proved. To prove the first statement, we construct a morphism from $\Lambda[q][h_1,...,h_{m-p}]$ to the dual presentation and show that its kernel is the claimed ideal of relations. 

{\it Acknowledgements:} I am indebted to S. Fomin and W. Fulton whose comments and suggestions enlightened the presentation of this paper. I am also thankful to A. L. Mare and A. Yong for some useful discussions and remarks.

\section{Factorial Schur functions}\label{fSchurs}
This section presents those properties of the factorial Schur
functions which are used later in the paper. The factorial Schur function $s_\lambda(x|a)$
is a homogeneous polynomial in two sets of variables:
$x=(x_1,...,x_p)$ and a doubly infinite sequence $a=(a_i)_{i \in
\z}$. An initial (non-homogeneous) version of these polynomials,
for $a_i = i-1$ if $1 \leqslant i \leqslant p$ and
$0$ otherwise, was first studied by L. Biedenharn and J. Louck in
\cite{BL}, then by W. Chen and J. Louck in \cite{CL}. The general
version was considered by I. Macdonald \cite{Mc2}, then studied
further in \cite{GG,Mo,MS} and \cite{Mc1}, Ch. I.3. These
functions, as well as a different version of them, called {\em
shifted Schur functions} (see \cite{O,OO}), play an important role
in the study of the center of the universal enveloping
algebra of $\frak{gl}(n)$, Capelli identities and quantum immanants. In a geometric
context, the factorial Schur functions appeared in \cite{KT1},
expressing the result of the localization of an equivariant
Schubert class to a $T-$fixed point.

For any variable $y$ and any sequence $(a_i)$ define the ``generalized factorial power": \[(y|a)^k=(y-a_1)\cdot \cdot \cdot
(y-a_k). \] Let $\lambda$ be a partition with at most $p$ parts.
Following \cite{MS}, define the factorial Schur function
$s_\lambda(x|a)$ to be
\[ s_{\lambda}(x|a)= \frac{\det[(x_j|a)^{\lambda_i+p-i}]_{1
\leqslant i,j \leqslant p}}{\det[(x_j|a)^{p-i}]_{1 \leqslant i,j \leqslant p}}. \]
Denote by $h_k(x|a)$ (respectively $e_k(x|a)$) the
factorial complete homogeneous Schur functions (resp. the
factorial elementary Schur functions). We adopt the usual
convention that $h_k(x|a)$ and $e_k(x|a)$ are equal to zero if $k$
is negative. For $k$ larger than $p$, $e_k(x|a)$ is set equal to
zero, also by convention. Let $\z[a]$ denote the ring of
polynomials in variables $a_i$ ($i$-integer). We start enumerating
the relevant properties of
the factorial Schur functions:\\

(A) (\textbf{Basis}) The factorial Schur functions $s_\lambda(x|a)$, where $\lambda$ has at most $p$ parts, form a $\z[a]-$basis for the ring of polynomials in $\z[a][x]$ symmetric in the
$x-$variables (\cite{Mc1}, I.3, pag. 55). This is a consequence
of the fact that
\[ s_\lambda(x|a) = s_\lambda(x) + \textrm{terms of lower degree in $x$.}
\] Then  \[ s_{\lambda}(x|a)s_\mu(x|a) = \sum_\nu
c_{\lambda\mu}^{\nu}(a) s_\nu(x|a) \] where the coefficients
$c_{\lambda\mu}^{\nu}(a)$ are homogeneous polynomials in $\z[a]$
of degree $|\lambda|+|\mu| - |\nu|$. They are equal to the usual 
Littlewood-Richardson coefficients $c_{\lambda\mu}^{\nu}$ if
$|\lambda|+|\mu| = |\nu|$ and are equal to $0$ if $|\lambda|+|\mu|
< |\nu|$ (cf. \cite{MS} \S 2).\\

(B) (\textbf{Vanishing Theorem}, \cite{MS}, Thm. 2.1, see also
\cite{O}) Let $\lambda=(\lambda_1,...,\lambda_p)$ and $\rho=(\rho_1,...,\rho_p)$ be two partitions of length at most $p$. Define the sequence $a_\rho$ by
$a_{\rho}=(a_{\rho_1+p},...,a_{\rho_p+1})$. Then

\begin{equation}\label{fvanishing} s_\lambda(a_\rho|a) = \left \{
\begin{array}{ll}
0 & \textrm{if } \lambda \subsetneq \rho \\
\prod_{(i,j) \in
\lambda}(a_{\lambda_i+p-i+1}-a_{p-\lambda_j^{'}+j}) & \textrm {if
} \lambda=\rho
\end{array} \right. \end{equation}
where $\lambda'$ is the conjugate partition of $\lambda$, and $(i,j) \in \lambda$ means that 
$j$ varies between $1$ and $\lambda_i$, if $\lambda_i >0$.
We recall a consequence of the Vanishing Theorem:
\begin{cor}\label{loc} The coefficients $c_{\lambda\mu}^\nu(a)$ satisfy the
following properties:

(i) $c_{\lambda\mu}^\nu(a) = 0$ if $\lambda$ or $\mu$ are not
included in $\nu$.

(ii) If the partitions $\mu$ and $\nu$ are equal, then
\[c_{\lambda\mu}^\mu(a) = s_\lambda(a_\mu |a). \]
\end{cor}
\begin{proof}(i) is part of the Theorem 3.1 in
\cite{MS} while (ii) follows from the proof of expression (10) in
\cite{MS}.\end{proof}

(C) (\textbf{Factorial Pieri-Chevalley rule}, see the proof of Prop. 3.2 in \cite{MS} or \cite{OO}, Thm. 9.1.) Let $(1)$ denote the
partition $(1,0,...,0)$ ($p$ parts) and let $\lambda$ be a
partition with at most $p$ parts. Then
\begin{equation}\label{fPieri} s_{(1)}(x|a)s_\lambda(x|a)=
\sum_{\mu \rightarrow \lambda} s_\mu (x|a) +
c_{(1),\lambda}^\lambda(a) s_\lambda(x|a) \end{equation} where
$\mu \rightarrow \lambda$ means that $\mu$ contains $\lambda$ and
has one more box than $\lambda$ (recall that, by definition, $\mu$
has at most $p $ parts). By Corollary \ref{loc},
$c_{(1),\lambda}^\lambda(a)=s_{(1)}(a_\lambda|a)$ and the last
expression turns out to be
\begin{equation}\label{boxcoeff} s_{(1)}(a_\lambda|a)= \sum_{i=1}^p
a_{\lambda_i+p+1-i} - \sum_{j=1}^p a_j.
\end{equation}

(D) (\textbf{Jacobi-Trudi identities}, see \cite{Mc1}, I.3, Ex.
20(c), pag. 56 or \cite{Mo}, Thm. 3.1.) Let $\tau^r a$ be the
sequence whose $n^{th}$ term is $a_{n+r}$ and let $\lambda$ be a
partition with at most $p$ parts. Then
\begin{eqnarray}\label{fjt} s_\lambda(x|a) & = & \det [
h_{\lambda_i-i+j}(x|\tau^{1-j}a)]_{1 \leqslant i,j \leqslant p} \\ & = & \det [
e_{\lambda_i^{'}-i+j}(x|\tau^{j-1}a)]_{1 \leqslant i,j \leqslant m-p},
\end{eqnarray}
where $\lambda'$ is the partition conjugate of $\lambda$.

The following proposition gives an inductive way of computing the
``shifted'' polynomials $h_k(x|\tau^sa)$, respectively
$e_k(x|\tau^sa)$, starting from the ``unshifted'' ones.
\begin{prop}\label{tauh} The following identities hold in $\z[a][x]$:
\begin{equation}\label{htau} h_{i+1}(x|\tau^{-1}a)=h_{i+1}(x|a) + (a_{i+p} - a_0)
h_i(x|a)
\end{equation}
\begin{equation}\label{etau} e_{j+1}(x|\tau a)=e_{j+1}(x|a)+ (a_1-a_{p-j+1})e_{j}(x|a).
\end{equation}
\end{prop}
\begin{proof} One uses the formulae
\begin{equation}\label{eq def of hk explicit} h_k(x|a) =
\sum_{1 \leqslant i_1 \leqslant ... \leqslant i_k \leqslant p}
(x_{i_1}-a_{i_1})(x_{i_2}-a_{i_2+1})\cdot ... \cdot
(x_{i_k}-a_{i_k+k-1}), \end{equation} 
\begin{equation}\label{eq def of ek expl} e_k(x|a) = \sum_{1 \leqslant i_1 < ... < i_k \leqslant p}
(x_{i_1}-a_{i_1})(x_{i_2}-a_{i_2 -1})\cdot ... \cdot
(x_{i_k}-a_{i_k+ k-1})
\end{equation}(cf. \cite{Mo}, eqs. (1.2) and (1.3)). Then the
computations are straightforward. \end{proof} By the Jacobi-Trudi formula, $h_i(x|\tau^{-s}a)$ (resp. $e_j(x|\tau^s a)$) generate the algebra of polynomials in $\z[a][x]$
symmetric in the $x-$variables. Then Prop. \ref{tauh} 
implies:
\begin{cor}\label{fgen} The factorial complete homogeneous (resp. elementary) symmetric functions
$h_i(x|a)$ for $1 \leqslant i$ (resp. $e_j(x|a)$ for $1 \leqslant
j \leqslant p$) generate the algebra of polynomials in $\z[a][x]$
symmetric in the $x-$variables.
\end{cor}

We will need the fact that the Jacobi-Trudi formula (2.5) generalizes to the case when $\lambda=(1)^k$ with $k>p$ (when, by convention, $e_k(x|a)=0$), and this is
the content of the next proposition (see also equation (6.10) in \cite{Mc2}):

\begin{prop}\label{eN} The following holds for any positive
integer $k>p$:
\[ \det\bigl(h_{1+j-i}(x|\tau^{1-j}a)\bigr)_{1 \leqslant i,j \leqslant k}=0. \]
\end{prop}
\begin{proof} Denote by $E_k(x|a)$ the determinant in question. Note that, if
$k \leqslant p$, this is equal to $e_k(x|a)$, by the
Jacobi-Trudi formula. We will need a formula proved in
\cite{Mc1}, I.3, pag.56, Ex. 20(b):
\begin{equation}\label{cauchy} \sum_{r=0}^p (-1)^r
e_r(x|a)h_{s-r}(x|\tau^{1-s}a)=0 \end{equation} for any positive
integer $s$. To prove the proposition, we use induction on $k \geqslant p+1$. Expanding $E_{p+1}(x|a)$
after the last column yields: \begin{eqnarray*} E_{p+1}(x|a)& = &
\sum_{r=0}^p (-1)^{r+1+p+1}h_{p+1-r}(x|\tau^{-p}a)e_r(x|a) \\ & =
& (-1)^{p+2}\sum_{r=0}^p (-1)^{r}h_{p+1-r}(x|\tau^{-p}a)e_r(x|a)\\
& = & 0 \end{eqnarray*} where the last equality follows from
(\ref{cauchy}) by taking $s=p+1$. Assume that $E_k(x|a)=0$ for all
$p<k<k_0$. Expanding the determinant defining
$E_{k_0}(x|t)$ after the last column, and using the induction
hypothesis, yields
\begin{eqnarray*} E_{k_0}(x|a)& = &
\sum_{r=0}^p (-1)^{r+1+k_0}h_{k_0-r}(x|\tau^{1-k_0}a)e_r(x|a) \\ &
=
& (-1)^{k_0+1}\sum_{r=0}^p (-1)^{r}h_{k_0-r}(x|\tau^{1-k_0}a)e_r(x|a)\\
& = & 0 \end{eqnarray*} using again (\ref{cauchy}) with $s=k_0$.
\end{proof}

\section{Equivariant quantum cohomology of the
Grassmannian}\label{sec eqq} In this section we recall some basic
properties of the equivariant quantum cohomology. As before, $X$
denotes the Grassmannian $Gr(p,m)$, and $\Lambda$ the polynomial
ring $\z[T_1,...,T_m]$. The ($T-$)equivariant quantum cohomology
of the Grassmannian, denoted by $\eqqhx$, is a deformation of both
equivariant and quantum cohomology rings (for details on the
latter cohomologies, see e.g. \cite{Mi1}). More precisely,
$\eqqhx$ is a graded, commutative, $\Lambda[q]-$algebra, where the degree of $q$
is equal to $m$, which has a $\Lambda[q]-$basis
$\{\sigma_\lambda\}$ indexed by the partitions $\lambda$ included in
the $p \times (m-p)$ rectangle. If
$\lambda=(\lambda_1,...,\lambda_p)$, the degree of
$\sigma_\lambda$ is equal to $|\lambda|=\lambda_1+...+\lambda_p$.
The multiplication of two basis elements $\sigma_\lambda$ and
$\sigma_\mu$ is given by the equivariant quantum
Littlewood-Richardson (EQLR) coefficients
$c_{\lambda,\mu}^{\nu,d}$, where $d$ is a nonnegative integer:
\[ \sigma_\lambda \circ \sigma_\mu = \sum_d \sum_\nu
q^d c_{\lambda,\mu}^{\nu,d} \sigma_\nu. \] Recall that the EQLR
coefficient $c_{\lambda,\mu}^{\nu,d}$ is a homogeneous polynomial
in $\Lambda$ of polynomial degree $|\lambda|+|\mu|-|\nu|-md$. If
$d=0$ one recovers the structure constant $c_{\lambda,\mu}^\nu$ in
the equivariant cohomology of $X$, and if the polynomial degree is
equal to $0$ (i.e. if $|\lambda|+|\mu|=|\nu|+md$) the EQLR
coefficient is equal to the ordinary $3-$pointed, genus $0$, Gromov-Witten
invariant $c_{\lambda,\mu}^{\nu,d}$. The latter is a nonnegative
integer equal to the number of rational curves in $X$ passing
through general translates of the Schubert varieties in $X$
corresponding to the partitions $\lambda,\mu$ and the dual of
$\nu$.

The geometric definition of these coefficients can be found in
\cite{Mi1}. In fact, for the purpose of this paper, the algebraic
characterization of the equivariant quantum cohomology from
Proposition \ref{cor char of eqq} below (which has a geometric
proof in {\it loc. cit.}), suffices. We only remark that the equivariant quantum Schubert
classes $\sigma_\lambda$ are determined by the equivariant Schubert classes $\sigma_\lambda^T$, determined in turn by the Schubert varieties
in $X$ defined with respect to the standard
flag.\begin{footnote}{Unlike the case of classical cohomology, in equivariant
cohomology the Schubert class determined
by a Schubert variety $\Omega_\lambda(F_\bullet)$, where
$F_\bullet$ is a fixed flag in $\cx^m$, depends on $F_\bullet$.}\end{footnote} The precise definition of $\sigma_\lambda^T$ is not presently needed, but it is given in \S \ref{sec eq giambelli}, where the equivariant cohomology ring is discussed in more detail. 

From now on we specialize the sequence $a=(a_i)_{i
\in \z}$ from the previous section to one, denoted $t=(t_i)_{i \in
\z}$, encoding the equivariant parameters $T_i$:
$$ t_i = \left \{
\begin{array}{ll}
T_{m-i+1} & \textrm{if } 1 \leqslant i \leqslant m \\
0 & \textrm { otherwise.}
\end{array} \right. $$
Using this sequence, we recall next the equivariant quantum
Pieri-Chvalley rule, as proved in \cite{Mi1}. Given a partition
$\lambda$, we denote by $\lambda^-$ the (uniquely determined)
partition obtained by removing $m-1$ boxes from the border rim of
$\lambda$ (recall that the border rim of a Young diagram is the
set of boxes that intersect the diagram's SE border - see also the
figure below). \vspace{1cm}
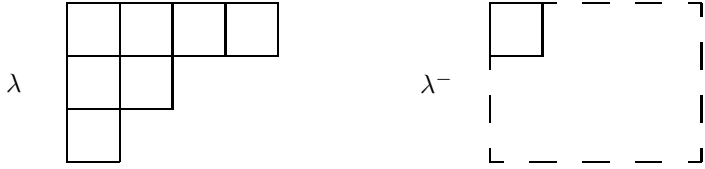
\begin{figure}[h!]
\begin{center}
\setlength{\unitlength}{2pt}
\begin{picture}(20,20)(0,0)

\thinlines

\put(-50,0){\line(0,1){30}} \put(-40,0){\line(0,1){30}}
\put(-30,10){\line(0,1){20}} \put(-20,20){\line(0,1){10}}
\put(-10,20){\line(0,1){10}} \put(-50,0){\line(1,0){10}}
\put(-50,10){\line(1,0){20}} \put(-50,20){\line(1,0){40}}
\put(-50,30){\line(1,0){40}}

\thicklines

\put(-60,15){\makebox(0,0){$\lambda$}} \thinlines
\put(30,20){\line(1,0){10}}

\put(30,30){\line(1,0){10}}

\put(30,20){\line(0,1){10}}

\put(40,20){\line(0,1){10}}

\put(20,15){\makebox(0,0){$\lambda^-$}}

\put(30,0){\dashbox{5}(40,30)}

\end{picture}
\end{center}\caption{Example: $p=3, m=7, \lambda=(4,2,1)$; $\lambda^-=(1)$.}
\end{figure}

\noindent If $\lambda=(\lambda_1,...,\lambda_p)$, note that
$\lambda^-$ exists only if $\lambda_1=m-p$ and $\lambda_p>0$.

\begin{prop}[Equivariant quantum Pieri-Chvalley rule - cf. \cite{Mi1}, Thm. 1]\label{eqq
Pieri} The following formula holds in $\eqqhx$:
 \[\sigma_\lambda \circ \sigma_{(1)} =
\sum_{\mu \rightarrow \lambda} \sigma_{\mu} + c_{\lambda,(1)
}^\lambda(t) \sigma_\lambda + q \sigma_{\lambda^-} \] where, by
the formula (\ref{boxcoeff}), $c_{\lambda,(1)}^\lambda(t)$ is
equal to
\[c_{\lambda,(1)}^\lambda(t)=\sum_{i=1}^p
T_{m-p+i-\lambda_i} - \sum_{j=m-p+1}^m T_j. \] The last term is
omitted if $\lambda^-$ does not exist.
\end{prop}

It turns out that the equivariant quantum Pieri-Chvalley rule
determines completely the equivariant quantum cohomology algebra,
in the following sense:

\begin{prop}[\cite{Mi1} Corollary 7.1]\label{cor char of eqq} Let $(A, \diamond)$
be a graded, commutative, associative $\Lambda[q]-$algebra with
unit such that:

\noindent 1. $A$ has an additive $\Lambda[q]-$basis
$\{s_\lambda\}$ (graded as usual).

\noindent 2. The equivariant quantum Pieri-Chevalley rule holds, i.e.
\[s_\lambda \diamond s_{(1)} =
\sum_{\mu \rightarrow \lambda} s_{\mu} + c_{\lambda,(1)
}^\lambda(t) s_\lambda + q s_{\lambda^-} \] where the last term is
omitted if $\lambda^-$ does not exist.

Then $A$ is canonically isomorphic to $\eqqh$, as
$\Lambda[q]-$algebras. \end{prop}

This proposition will be the main tool in proving the presentation
and equivariant quantum Giambelli formula from the next section.

\section{Proof of the Theorem}\label{sec proof}
The strategy for the proof is to ``guess'' candidates for the
presentation and for the polynomial representatives, using the insight provided by
the similar results in quantum cohomology (see e.g. \cite{BCF}) and some related results in equivariant cohomology (\cite{KT1} \S 6). Then one attempts to prove that the guessed polynomials form a $\Lambda[q]-$basis in the
candidate presentation, and they multiply according to the EQ
Pieri-Chevalley rule. Proposition \ref{cor char of eqq} will
ensure that the guessed algebra will be canonically
isomorphic to $\eqqhx$ and that the polynomials considered will
represent the equivariant quantum Schubert classes.  

It turns out that each of the quantum presentations from
\cite{BCF} (the usual one, involving the $h$ variables, and the
``dual'' one, involving the variables $e$) implies an equivariant
quantum presentation (see respectively Theorems \ref{eqq} and
\ref{eqq-el} below). The equivariant generalizations are obtained
by taking the factorial versions, via the factorial Jacobi-Trudi
formula (\S \ref{fSchurs}, property D) of all the expressions
involved in the original quantum presentations.

Before stating the first result, we recall the notation from the
introduction: $h_1,...,h_{m-p}$ and $e_1,...,e_p$ denote two sets
of indeterminates; the definitions of $\tau^{-s}h_i$, $\tau^se_j$
and of $H_k$ ($m-p < k \leqslant m$) respectively $E_k$ ($p < k
\leqslant m$) are those given in the equations (\ref{eq def htau
introx}),(\ref{eq def etau introx}) and (\ref{eq def Hk introx}),(\ref{eq def Ek introx}) above. For $\lambda$ in the $p
\times (m-p)$ rectangle recall that:
\begin{equation}\label{sh} s_\lambda=
\det(\tau^{1-j}h_{\lambda_i+j-i})_{1 \leqslant i,j \leqslant p}
\end{equation} respectively
\begin{equation}\label{se} \widetilde{s}_\lambda=\det(\tau^{j-1}e_{\lambda_i'+j-i})_{1
\leqslant i,j \leqslant m-p}, \end{equation} (cf. (\ref{sh introx})
and (\ref{se introx})) with the usual conventions that $h_k=0$ for
$k<0$ and $k>m-p$ respectively $e_i=0$ if $i<0$ or $i>p$.
Before proving the theorem, we need a Nakayama-type result, which
will be used several times in the paper:
\begin{lemma}[cf. \cite{E}, Exerc. 4.6]\label{Nak}
Let $M$ be an $R-$algebra graded by nonnegative integers. Assume
that $R$ is also graded (by nonnegative integers) and let $I$ be a
homogeneous ideal in $R$ consisting of elements of positive
degree. Let $m_1,...,m_k$ be homogeneous elements whose images
generate $M/IM$ as an $R/I$-module. Then $m_1,...,m_k$ generate
$M$ as an $R-$module.
\end{lemma}
\begin{proof} Let $m$ be a nonzero homogeneous element of $M$. We use induction
on its degree. Assume $\deg m = 0$. The hypothesis implies that
\begin{equation}\label{m} m = r_1 m_1 + ... + r_k m_k \mod{IM}  \end{equation}
where $r_i$ are elements in $R$. Since $I$ contains only elements
of positive degree, it follows that the equality holds in $M$ as
well. Let now $\deg m > 0$. Writing $m$ as in (\ref{m}), implies
that \[m - \sum_i r_i m_i = \sum_j a_j m'_j \] for some (finitely
many) $a_j \in I$ and $m'_j \in M$. Again, since $I$ contains only
elements of positive degree, $\deg m'_j < \deg m$ for each $j$.
The induction hypothesis implies that each $m'_j$ is an
$R-$combination of $m_i$'s, which finishes the proof.
\end{proof}
We prove next the ``dual version" statement from the main theorem.

\begin{thm}\label{eqq-el} There exists a canonical isomorphism of $\Lambda[q]-$algebras
\[ \Lambda[q][e_1,...,e_p]/ \langle H_{m-p+1},...,H_m+(-1)^p q\rangle
\longrightarrow \eqqhx, \] sending $e_i$ to $\sigma_{(1)^i}$ and
$\widetilde{s}_\lambda$ to the equivariant quantum
Schubert class $\sigma_\lambda$.
\end{thm}
\begin{proof} Note first that
\[ e_j(x|t)= e_j(x) + f(t,x), \] where $f(t,x)$ is a homogeneous
polynomial in the variables $x$ and $t$, but of degree in the
variables $x$ less than $j$. Since the usual elementary symmetric
functions $e_1(x),...,e_p(x)$ are algebraically independent over
$\z$, it follows that the elementary factorial Schur functions
$e_1(x|t),...,e_p(x|t)$ are algebraically independent over
$\Lambda$. Then there is a canonical isomorphism
\[ \Lambda[q][e_1(x|t),...,e_p(x|t)] \to \Lambda[q][e_1,...,e_p], \]
sending $s_\lambda(x|t)$ to $\widetilde{s}_{\lambda}$, and $h_{m-p+i}(x|t)$ ($1 \leqslant i \leqslant p$) to $H_{m-p+i}$, by the
factorial Jacobi-Trudi identity. This induces an isomorphism
between
\[A:=\Lambda[q][e_1(x|t),...,e_p(x|t)]/ \langle h_{m-p+1}(x|t),...,
h_m(x|t)+(-1)^p q\rangle\]
and
\[\Lambda[q][e_1,...,e_p]/ \langle H_{m-p+1},...,H_m+(-1)^p
q\rangle.\] By Prop. \ref{cor char of eqq}, it remains to show
that the images of $s_\lambda(x|t)$ in $A$, as
$\lambda$ varies over the partitions included in the $p \times
(m-p)$ rectangle, form a $\Lambda[q]$-basis of $A$, satisfying the
equivariant quantum Pieri-Chevalley rule.

{\it Generating set.} This follows from Lemma \ref{Nak}, applied
to $M=A$, $R=\Lambda[q]$ and $I$ the ideal generated by $q$ and
$T_1,...,T_p$ (in which case $M/I$ is the classical cohomology of $X$).

{\it Linear independence.} Assume that $\sum q^{d_\lambda}
c_\lambda s_\lambda(x|t) = 0$ in $A$, for $c_\lambda$ in $\Lambda$, where $\lambda$ is included in the $p \times (m-p)$ rectangle.
This implies that $\sum q^{d_\lambda} c_\lambda s_\lambda(x|t)$ is
in the ideal generated by $h_{m-p+1}(x|t),...,h_m(x|t)+(-1)^p q$.
By Cor. \ref{loc} (1), any element of this ideal can be written as:
\[ \sum_\mu q^{d'_\mu}c'_\mu s_\mu(x|t)+\sum_\nu q^{d_\nu'}c_\nu'' s_\nu(x|t)(h_m(x|t)+(-1)^pq) \] where $\mu,\nu$  have at most $p$ parts, $\mu$ is {\it outside} the $p \times (m-p)$ rectangle, and $c'_\mu,c''_\nu$ are in $\Lambda$. Note that $s_\nu(x|t) h_m(x|t)$ expands also into a sum of factorial Schur functions indexed by partitions outside the $p \times (m-p)$ rectangle. Since the factorial Schur functions form a $\Lambda[q]-$basis for the polynomials in
$\Lambda[q][x_1,...,x_p]$ symmetric in the $x-$variables, it
follows that all $c_\lambda$ (and $c'_\mu,c''_\nu$) must be equal to zero,
as desired.

{\it Equivariant quantum Pieri-Chevalley.} The factorial
Pieri-Chevalley rule (\S \ref{fSchurs}, Property (C)) states that
if $\lambda$ is included in the $p \times (m-p)$ rectangle, then
\[ s_\lambda(x|t) \cdot s_{(1)}(x|t) = \sum_{\mu}
s_\mu(x|t) + c_{(1),\lambda}^\lambda(t) s_\lambda(x|t) +
s_{\overline{\lambda}}(x|t)
\] where $\mu$ runs over all partitions in the $p \times (m-p)$
obtained from $\lambda$ by adding one box; the last term is
omitted if $\lambda_1$, the first part of $\lambda$, is not equal
to $m-p$. If $\lambda_1=m-p$ then
$\overline{\lambda}=(\overline{\lambda}_1,...,\overline{\lambda}_p)$
is given by adding a box to the first row of $\lambda$, i.e.
$\overline{\lambda}_1 = m-p+1$ and $\overline{\lambda}_i =
\lambda_i $ for $i \geqslant 2$. Since the images of $s_\lambda(x|t)$, as
$\lambda$ varies in the $p \times (m-p)$ rectangle, form a
$\Lambda[q]-$basis for $A$, it is enough to show that
\[s_{\overline{\lambda}}(x|t) = q s_{\lambda^-}(x|t)
\mod{J} \] where $J$ is the ideal generated by
$h_{m-p+1}(x|t),...,h_m(x|t)+(-1)^p q$. %Extend the definition of
%$H^T_i$ for $i \leqslant m-p$.
By the factorial Jacobi-Trudi formula (\S \ref{fSchurs}, Property
(D)), it follows that
\begin{equation}\label{sdet} s_{\overline{\lambda}}(x|t)= \det \begin{pmatrix}
h_{m-p+1}(x|t) &
h_{m-p+2}(x|\tau^{-1}t) & ... & h_{m}(x|\tau^{1-p}t) \\
h_{\lambda_2-1}(x|t) & h_{\lambda_2}(x|\tau^{-1}t) & \hdots &
h_{\lambda_2+ p-2}(x|\tau^{1-p}t)
\\ \vdots & \vdots & \vdots & \vdots \\ \hdots & \hdots & \hdots
& h_{\lambda_p}(x|\tau^{1-p}t)
\end{pmatrix}. \end{equation}

We analyze next the first row of this determinant. \\ %Equation (\ref{htau}) implies %that \begin{equation}\label{Htau}
%h_j(x|\tau^{-s}t)= h_j(x|\tau^{-s+1}t) + (t_{j+p-s}-t_{-s+1})
%h_{j-1}(x|\tau^{-s+1}t).\end{equation}

\noindent {\bf Claim.} Let $i,j$ be two integers such that
$2 \leqslant j \leqslant i \leqslant p$. Then
\begin{equation*}\label{eq equality}
h_{m-p+i}(x|\tau^{1-j}t)=h_{m-p+i}(x|\tau^{1-(j-1)}t).\end{equation*}

\begin{proof} [Proof of the Claim] Equation (\ref{htau}) implies that \begin{equation*}\label{Htau}
h_j(x|\tau^{-s}t)= h_j(x|\tau^{-s+1}t) + (t_{j+p-s}-t_{-s+1})h_{j-1}(x|\tau^{-s+1}t),\end{equation*}
hence,
\[ h_{m-p+i}(x|\tau^{1-j}t)=
h_{m-p+i}(x|\tau^{1-(j-1)}t)+(t_{m+1+i-j}-t_{2-j})h_{m-p+i-1}(x|\tau^{1-(j-1)}t).\]
The Claim follows then from the definition of $(t_i)$, since
$t_{m+1+i-j}=t_{2-j}=0$. \end{proof}

It follows that for any integer $1 \leqslant i \leqslant
p$,
\begin{equation}\label{eq equality dec} h_{m-p+i}(x|\tau^{1-i}t) =
h_{m-p+i}(x|t). \end{equation} In particular, \[
h_{m-p+i}(x|\tau^{1-i}t) = 0 \mod{J}
\] if $1 \leqslant i \leqslant p-1$, and \[h_m(x|\tau^{1-p}t) = (-1)^{p+1}q \mod{J}.\]
Therefore,
expanding the determinant in (\ref{sdet}) after the first row,
yields:
\begin{equation}\label{slambda-} s_{\overline{\lambda}}(x|t)= (-1)^{p+1} (-1)^{p+1} q
\det \begin{pmatrix} h_{\lambda_2-1}(x|t) &
h_{\lambda_2}(x|\tau^{-1}t) & \hdots & h_{\lambda_2+
p-3}(x|\tau^{2-p}t)
\\ \vdots & \vdots & \vdots & \vdots \\ \hdots & \hdots & \hdots
& h_{\lambda_{p}-1}(x|\tau^{2-p}t)
\end{pmatrix} \end{equation}
in $A$. If $\lambda_p=0$, the last row of the determinant in
(\ref{slambda-}) contains only zeroes; if $\lambda_p > 0 $, the
determinant is equal to $s_{\lambda^-}(x|t)$, by the Jacobi-Trudi
formula. Summarizing, $s_{\overline{\lambda}}(x|t)$ is equal to $q
s_{\lambda^-}(x|t)$ in $A$, or it is equal to zero if $\lambda^-$
does not exists. This finishes the proof of the equivariant
quantum Pieri-Chevalley rule, hence also the proof of the
theorem.\end{proof}

We are ready to prove the first part of the main result, which
involves the $h$ variables. We use the notation preceding Thm.
\ref{eqq-el} above.

\begin{thm}\label{eqq} There exist a canonical isomorphism of $\Lambda[q]-$algebras
\[ \Lambda[q][h_1,...,h_{m-p}]/ \langle E_{p+1},....,E_m +(-1)^{m-p}q \rangle
\longrightarrow \eqqhx, \] such that $h_j$ is sent to
$\sigma_{(j)}$ and $s_\lambda$ to the equivariant quantum Schubert
class $\sigma_\lambda$. \end{thm}
\begin{proof} Consider the $\Lambda[q]-$algebra morphism \[ \Psi:\Lambda[q][h_1,...,h_{m-p}] \to
\Lambda[q][e_1(x|t),...,e_p(x|t)]/ \langle h_{m-p+1}(x|t),...,
h_m(x|t)+(-1)^p q\rangle \] sending $h_k$ to the image of
$h_k(x|t)=\det(e_{1+j-i}(x|\tau^{1-j}t))_{1 \leqslant i,j \leqslant k}$. Recall that the last quotient is denoted by $A$ and it
is canonically isomorphic to $\eqqhx$, by the previous proof. We
will show that the images under $\Psi$ of
$E_{p+1},E_{p+2},...,E_{m-1},E_m+(-1)^{m-p}q$ are equal to zero in
$A$ (where $E_i$ is defined by equation (\ref{eq def Ek
introx})). First, we need the following claim:\\

\noindent {\bf Claim.} The following formulae hold in $A$:
\begin{equation}\label{eq psi<m}\Psi(\tau^{-s}h_j)= h_j(x|\tau^{-s}t), \end{equation}
for any nonnegative integers $s$ and $j$ with $ j < m$, and
\begin{equation}\label{eq psi=m}\Psi(\tau^{-(m-1)}h_{m})=
h_m(x|\tau^{-(m-1)}t) + (-1)^p q .\end{equation}
\begin{proof}[Proof of the Claim] By definition, both $\tau^{-s}h_j$ and $h_j(x|\tau^{-s}t)$ satisfy the same recurrence relations (given respectively by the equations (\ref{eq def htau introx}) and (\ref{htau})). This implies that there exist polynomials $P_1(t),...,P_s(t)$ in $\Lambda$, with $\deg P_k(t)=k$, such that 
\[ \tau^{-s}h_j=h_j + \sum_{k=1}^s P_k(t) h_{j-k}, \] respectively
\begin{equation}\label{eq expansion}h_j(x|\tau^{-s}t)=h_j(x|t) + \sum_{k=1}^s P_k(t) h_{j-k}(x|t). \end{equation}
If $j \leqslant m-p$, then $\Psi(h_j)=h_j(x|t)$ in $A$, by the definition of $\Psi$, thus \begin{equation}\label{eq star}\Psi(\tau^{-s}h_j)= h_j(x|\tau^{-s}t). \end{equation} If $m-p+1 \leqslant j < m$, $h_j=0$ by convention, whereas $h_j(x|t)=0$ in $A$, so equation (\ref{eq star}) also holds in this case. If $s=m-1$ and $j=m$, we have \begin{equation*}\begin{split}  \Psi(\tau^{-(m-1)}h_m) & =  \Psi(h_m + \sum_{k=1}^{m-1} P_k(t) h_{m-k}) \\ & = \Psi(h_m) + \sum_{k=1}^{m-1}P_k(t) \Psi(h_{m-k})\\ & = h_m(x|t)+(-1)^p q + \sum_{k=1}^{m-1}P_k(t) h_{m-k}(x|t) \\ & = h_m(x|\tau^{-(m-1)}t)+(-1)^pq,\end{split}\end{equation*}
where the third equality follows from the fact that $h_m=0$ and $h_m(x|t)+(-1)^pq=0$ in $A$; the fourth equality follows from the expansion (\ref{eq expansion}) of $h_m(x|\tau^{-(m-1)}t)$.\end{proof}
By definition, $\Psi(E_i)$ is equal to the image in $A$, through $\Psi$, of
\[ \det \begin{pmatrix} h_1 & \tau^{-1}h_2 & ...& \tau^{-(s-1)}h_{s} & ... &
\tau^{-(i-1)}h_{i} \\1 & \tau^{-1}h_1 & ... & ...& ... &
\tau^{-(i-1)}h_{i-1} \\ \vdots & \vdots & \vdots & \vdots & \vdots
& \vdots \\ 0 & ... & ... & 0 & 1 & \tau^{-(i-1)}h_1
\end{pmatrix}. \]
If $p+1 \leqslant i <m$, this determinant contains only
$\tau^{-s}h_j$ with $j<m$. Then, by the equation (\ref{eq psi<m}) from the claim,  
$\Psi(E_i)$ is the image in $A$ of the determinant
\[ \det \begin{pmatrix} h_1(x|t) & h_2(x|\tau^{-1}t) & ...& h_{s}(x|\tau^{-(s-1)}t) & ... &
h_{i}(x|\tau^{-(i-1)}t) \\1 & h_1(x|\tau^{-1}t) & ... & ...& ... &
h_{i-1}(x|\tau^{-(i-1)}t) \\ \vdots & \vdots & \vdots & \vdots &
\vdots & \vdots \\ 0 & ... & ... & 0 & 1 & h_1(x|\tau^{-(i-1)}t)
\end{pmatrix} \] which, by Proposition \ref{eN}, is equal to zero, since $i >p$.
To compute the image of $\Psi(E_m)$ we use both the equations (\ref{eq  psi<m}) and (\ref{eq psi=m}). Then $\Psi(E_m)$ is equal to the image in $A$ of \\
$$\begin{array}{c}
  \det \begin{pmatrix} h_1(x|t) & h_2(x|\tau^{-1}t) & ...&
h_{s}(x|\tau^{-(s-1)}t) & ... & h_{m}(x|\tau^{-(m-1)}t) \\ \vdots
& \vdots & \vdots & \vdots & \vdots & \vdots \\ 0 & ... & ... & 0
& 1 & h_1(x|\tau^{-(m-1)}t)
\end{pmatrix} + \\
 \det \begin{pmatrix} 0& 0 & ... & 0 & (-1)^p q \\ 1 &
h_1(x|\tau^{-1}t) & ... & ... & h_{m-1}(x|\tau^{-(m-1)}t) \\
\vdots & \vdots & \vdots & \vdots & \vdots \\ 0 & ... & 0 & 1 &
h_1(x|\tau^{-(m-1)}t)
\\ \end{pmatrix}  \\
\end{array}$$
The first determinant is equal to zero, by Prop. \ref{eN}, and the second is $(-1)^{m+1}\cdot (-1)^p
q$. It follows that $\Psi(E_m)+(-1)^{m-p} q$ is equal to zero in
$A$, as claimed. Thus $\Psi$ induces a $\Lambda[q]-$algebra morphism
\[ \Psi': \Lambda[q][h_1,...,h_{m-p}]/ \langle E_{p+1},....,E_m +(-1)^{m-p}q \rangle
\longrightarrow A . \] Note that $\tau^{-s}h_m$ does not appear in
the determinant defining $s_\lambda$, for any $\lambda$ in the $p
\times (m-p)$ rectangle, so $\Psi'$ sends $s_\lambda$ to the image
of $s_{\lambda}(x|t)$ in $A$. Applying Lemma \ref{Nak} with $M=
\Lambda[q][h_1,...,h_{m-p}]/ \langle E_{p+1},....,E_m +(-1)^{m-p}q
\rangle$ and $I$ the ideal generated by $q$ and $T_1,...,T_m$
implies that the polynomials $s_\lambda$ generate
$\Lambda[q][h_1,...,h_{m-p}]/ \langle E_{p+1},....,E_m
+(-1)^{m-p}q \rangle$ as a $\Lambda[q]$-module. Since their images
through $\Psi'$ form a $\Lambda[q]-$basis, they must form a
$\Lambda[q]-$basis, as well. Hence $\Psi'$ is an isomorphism, as
desired. \end{proof}

{\it Remark:} The Theorems \ref{eqq-el} and \ref{eqq} are proved without using the corresponding results from quantum cohomology. In particular, we obtain a new proof for Bertram's quantum Giambelli formula (see \cite{Be}). (Recall that the quantum cohomology ring of $X$ is a graded $\z[q]-$algebra isomorphic to $\eqqhx/(T_1,...,T_m)$, hence the quantum Giambelli formula is obtained by taking $T_1=...=T_m=0$ in the determinants from the above-mentioned theorems.)  

\section{Giambelli formulae in equivariant cohomology}\label{sec eq giambelli}
The goal of this section is to state the equivariant Giambelli formulae implied by their equivariant quantum counterparts from the previous section. We will also use this opportunity to define rigourously the equivariant Schubert classes involved, and provide, without proof, a geometric interpretation for the factorial Schur functions.  

Let $T$ be the usual torus, and let $ET \to BT$ be the universal $T-$bundle. If $X$ is a topological space with a $T-$action, there is an induced $T-$action on $ET \times X$ given by $t\cdot (e,x)=(et^{-1},tx)$. The (topological) quotient space $(ET \times X)/T$ is denoted by $X_T$. By definition, the ($T-$)equivariant cohomology of $X$, denoted $H^*_T(X)$, is equal to the usual cohomology of $X_T$. The $X-$bundle $X_T \to BT$ gives $H^*_T(X)$ the structure of a $\Lambda-$algebra, where $\Lambda$ denotes the equivariant cohomology of a point $H^*_T(pt)=H^*(BT)$.

Let now $X$ be the Grassmannian of subspaces of dimension $p$ in $\cx^m$ with the $T-$action induced from the usual $GL(m)-$action. We define next the equivariant Schubert classes which determine the equivariant quantum classes $\sigma_\lambda$ used in previous sections (see also \cite{Mi1}). Let  \[ F_\bullet: (0) \subset F_1 \subset
... \subset F_m = \cx^m \] be the standard flag, so $F_i = \langle f_1,...,f_i
\rangle$ and $f_i=(0,...,1,...,0)$ (with $1$ in the $i-$th position). 
If $\lambda=(\lambda_1,...,\lambda_p)$ is a partition included in
the $p \times (m-p)$ rectangle, define the Schubert variety
$\Omega_\lambda(F_\bullet)$ by
\begin{equation}\label{schvdef} \Omega_\lambda(F_\bullet) = \{ V \in
X : \dim V \cap F_{m-p+i-\lambda_i} \geqslant i \}.
\end{equation} Since the flag $F_\bullet$ is $T$-invariant, the
Schubert variety $\Omega_\lambda(F_\bullet)$ will be $T-$invariant
as well, so it determines a Schubert class
$\sigma^T_\lambda$ in $H^{2|\lambda|}_T(X)$.
The following result is a consequence of the Theorems
\ref{eqq-el} and \ref{eqq} (the notation is 
from the
previous section):

\begin{cor}\label{preq} (a)  There exists a canonical isomorphism of $\Lambda-$algebras
\[ \Lambda[h_1,...,h_{m-p}]/ \langle E_{p+1},....,E_m \rangle
\longrightarrow H^*_T(X), \] sending $h_j$ to $\sigma_{(j)}^T$ and $s_\lambda$
is the equivariant Schubert class $\sigma_\lambda^T$.

(b) There exists a canonical isomorphism of $\Lambda-$algebras
\[ \Lambda[e_1,...,e_p]/ \langle H_{m-p+1},...,H_m\rangle
\longrightarrow H^*_T(X), \] sending $e_i$ to $\sigma_{(1)^i}^T$ and
$\widetilde{s}_\lambda$ to the equivariant Schubert
class $\sigma_\lambda^T$.
\end{cor}
\begin{proof} It is known (see e.g. \cite{Mi1}) that there is a canonical
isomorphism of $\Lambda-$algebras
\[ \eqqhx/\langle q \rangle \longrightarrow H^*_T(X) \]
sending the equivariant quantum Schubert class $\sigma_\lambda$
from the previous section to $\sigma^T_\lambda$. Then the
Corollary follows from the Theorems \ref{eqq-el} and \ref{eqq}.
\end{proof}

{\it Remarks:} 1. The proof of the Corollary can be given without
using the equivariant quantum cohomology. There is an analogue of Prop. \ref{cor char of eqq}, stating that the Pieri-Chevalley rule determines the equivariant cohomology algebra.  Then a ``strictly equivariant" proof of Cor. \ref{preq} can be obtained by taking $q=0$ in all the assertions
from the previous section.

2. The fact that the factorial Schur functions represent the equivariant Schubert classes can be
also be deduced, indirectly, by combining the fact that the double Schubert polynomials represent the equivariant Schubert classes in the complete flag variety (cf. \cite{Bi} and \cite{A}) and that, when indexed by a Grassmannian permutation, these polynomials are actually factorial Schur functions. The latter holds because the vanishing property characterizing the factorial Schur functions (\S \ref{fSchurs}, property (B)), is also satisfied by the double Schubert polynomials in question (see \cite{L}, pag. 33). However, the details of this connection are missing from the literature.

3. It is well known that the equivariant Schubert classes are determined by their restriction to the torus fixed points in $X$. Formulae for such restrictions have been obtained by A. Knutson - T. Tao in \cite{KT1} and, recently, by V. Lakshmibai - K.N. Raghavan - P. Sankaran in \cite{LRS}.

\subsection{A geometric interpretation of the factorial Schur functions}
Consider the tautological short exact sequence on $X$:
\begin{equation}\label{tautseq} 0 \longrightarrow S \longrightarrow
V \longrightarrow Q \longrightarrow 0, \end{equation} which is clearly
$T-$equivariant. Let $-x_1,...,-x_p$ be the {\em equivariant} Chern roots of the bundle $S$.
There is a weight space decomposition of the trivial (but not equivariantly trivial) vector bundle $V$ into a sum of $T-$equivariant line bundles:
\[ V = L_1 \oplus ... \oplus L_m. \] Let $-T_i$ be the equivariant first Chern class of $L_i$.\begin{footnote}{All the minus signs are for positivity reasons. It turns out, for example, that $c_1^T(L_i)$ is the Chern class of $\mathcal{O}_{\mathbb{P}^\infty}(-1)$ (see e.g. \cite{Mi2} \S 7).}\end{footnote} Define the sequence $(t_i)$ as usual, using the formula from \S \ref{subsec statement of res}.

\begin{prop} \label{geomschur} In $H^*_T(X)$, the equivariant Schubert class $\sigma_\lambda^T$ is equal to the factorial Schur polynomial $s_\lambda(x|t)$. \end{prop}
\begin{proof}[Idea of proof] The equivariant Schubert class $\sigma_\lambda^T$ is a cohomology class on the infinite dimensional space $X_T$. The first step of the proof is to approximate this class by a class $(\sigma_\lambda)_{T,n}$ on a finite-dimensional ``approximation" $X_{T,n}$ ($n \gg 0$) of $X_T$. This is standard (see e.g. \cite{Br2} or \cite{Mi1}) and uses the $T-$bundle $(\cx^{n+1}\setminus 0)^m \to (\mathbb{P}^n)^m$ which approximates the universal $T-$bundle $ET \to BT$. Then $X_{T,n}:=(ET_n \times X)/T$, which, in fact, is equal to the Grassmann bundle $\mathbb{G}\bigl(p,\mathcal{O}_{(1)}(-1) \oplus ... \oplus \mathcal{O}_{(m)}(-1)\bigr)$, where $\mathcal{O}_{(i)}(-1)$ denotes the tautological line bundle over the $i-$th component of $(\mathbb{P}^n)^m$. Using this procedure one obtains the class $(\sigma_\lambda)_{T,n}$ as the cohomology class determined by the subvariety $(\Omega_\lambda)_{T,n}$ of $X_{T,n}$. The second step is to use the definition of $(\Omega_\lambda)_{T,n}$ to realize it as the degeneracy locus from \cite{F3}, Thm. 14.3, whose cohomology class is given as a certain determinantal formula in the Chern classes of the vector bundles $S_{T,n},V_{T,n}$ and $Q_{T,n}$ on $X_{T,n}$ induced by the tautological sequence on $X$. Finally, one proves that the determinant in question is equal to the claimed factorial Schur polynomial, which ends the proof.\end{proof}

\bibliographystyle{plain}
\bibliography{bibliography}

\end{document}